\documentclass[12pt]{article}
\usepackage{amssymb,amsfonts,amsmath, psfrag,eepic,colordvi,graphicx,epsfig,ytableau}
\usepackage{amssymb,latexsym,graphics}
\usepackage{MnSymbol}
\usepackage[enableskew]{youngtab}
\usepackage{float}
\usepackage{tikz}
\topskip  
\newcommand{\rmnum}[1]{\romannumeral #1}

 \numberwithin{equation}{section}
\newtheorem{theorem}{Theorem}[section]

\newtheorem{corollary}[theorem]{Corollary}

\newtheorem{conjecture}[theorem]{Conjecture}

\newtheorem{lemma}[theorem]{Lemma}

\begin{document}
\parskip 6pt

\pagenumbering{arabic}
\def\sof{\hfill\rule{2mm}{2mm}}
\def\ls{\leq}
\def\gs{\geq}
\def\SS{\mathcal S}
\def\qq{{\bold q}}
\def\MM{\mathcal M}
\def\TT{\mathcal T}
\def\EE{\mathcal E}
\def\lsp{\mbox{lsp}}
\def\rsp{\mbox{rsp}}
\def\pf{\noindent {\it Proof.} }
\def\mp{\mbox{pyramid}}
\def\mb{\mbox{block}}
\def\mc{\mbox{cross}}
\def\qed{\hfill \rule{4pt}{7pt}}
\def\block{\hfill \rule{5pt}{5pt}}

\begin{center}
{\Large\bf  Quasi-Stirling Polynomials on Multisets  }
\end{center}

\begin{center}
{\small  Sherry H.F.  Yan, Xue Zhu  }

 Department of Mathematics\\
Zhejiang Normal University\\
 Jinhua 321004, P.R. China

 hfy@zjnu.cn

\end{center}

\noindent {\bf Abstract.} A permutation $\pi$ of a multiset is said to be a {\em quasi-Stirling
} permutation  if  there does not exist four indices $i<j<k<\ell$ such that $\pi_i=\pi_k$ and $\pi_j=\pi_{\ell}$. For a multiset $\mathcal{M}$, denote by   $\overline{\mathcal{Q}}_{\mathcal{M}}$ the set of quasi-Stirling permutations of $\mathcal{M}$.  The   {\em  qusi-Stirling polynomial}   on the multiset $\mathcal{M}$ is defined by  $  \overline{Q}_{\mathcal{M}}(t)=\sum_{\pi\in \overline{\mathcal{Q}}_{\mathcal{M}}}t^{des(\pi)}$, where $des(\pi)$ denotes the number of descents of $\pi$.  By employing generating function arguments,  Elizalde  derived an   elegant identity involving   quasi-Stirling polynomials on the multiset $\{1^2, 2^2, \ldots, n^2\}$,   in  analogy to the identity on   Stirling polynomials.  In this paper,   we     derive an identity  involving quasi-Stirling polynomials  $\overline{Q}_{\mathcal{M}}(t)$ for any multiset $\mathcal{M}$, which is  a generalization of  the identity on Eulerian polynomial and Elizalde's identity on quasi-Stirling polynomials on the multiset $\{1^2, 2^2, \ldots, n^2\}$.    We provide a combinatorial proof  the identity  in terms of certain ordered labeled trees. Specializing  $\mathcal{M}=\{1^2, 2^2, \ldots, n^2\}$    implies a combinatorial proof of Elizalde's identity   in answer to the problem posed by Elizalde. As an application, our identity enables us to show that   the quasi-Stirling polynomial  $\overline{Q}_{\mathcal{M}}(t)$   has only real roots and  the coefficients of $\overline{Q}_{\mathcal{M}}(t)$  are unimodal  and log-concave    for any multiset $\mathcal{M}$,  in   analogy to
 Brenti's result  for Stirling polynomials on   multisets.




\section{Introduction}

Let $\mathcal{M}=\{1^{k_1}, 2^{k_2}, \ldots, n^{k_n}\}$ be a multiset where $k_i$ is the number of occurrences of $i$ in  $\mathcal{M}$ and $k_i\geq 1$.
  A permutation $\pi$ of a multiset $\mathcal{M}$ is said to be  a {\em Stirling } permutation if   $i<j<k$ and $\pi_i=\pi_k$, then $\pi_j>\pi_i$.
  Stirling permutations were originally introduced by Gessel and Stanley \cite{Gessel} in the case of the multiset $\mathcal{M}=\{1^{2}, 2^{2}, \ldots, n^{2}\}$. Denote by $\mathcal{Q}_{\mathcal{M}}$ the set of Stirling permutations of $\mathcal{M}$.

 For a permutation $\pi=\pi_1\pi_2\ldots \pi_n$,   the number of descents of $\pi$, denoted by  $des(\pi)$, is defined by
 $$
 des(\pi)=|\{i\mid \pi_i>\pi_{i+1}, 1\leq i\leq n-1\}|+1.
 $$
  The polynomial  $$ Q_{\mathcal{M}}(t)=\sum_{\pi\in \mathcal{Q}_{\mathcal{M}}}t^{des(\pi)}$$ is called
the  {\em Stirling polynomial } on the multiset $\mathcal{M}$. When $\mathcal{M}=\{1,2,\ldots, n\}$, the Stirling polynomial reduces to Eulerian polynomial $A_n(t)$.
Let
$$
\sum_{m=0}^{\infty}B_{\mathcal{M}}(m)t^m={Q_{\mathcal{M}}(t)\over (1-t)^{K+1}},
$$
where $\mathcal{M}=\{1^{k_1}, 2^{k_2}, \ldots, n^{k_n}\}$ and  $K=k_1+k_2+\ldots+k_n$.
Gessel and Stanley \cite{Gessel}
proved that  $B_{\mathcal{M}}(m)=S(n+m,m)$  when  $\mathcal{M}=\{1^2, 2^2, \ldots, n^2\}$,
where  $S(n+m,m)$ is the number of partitions of $[m+n]$ into $m$ blocks.
Brenti \cite{Brenti1, Brenti2} investigated  Stirling permutations   of an arbitrary  multiset $\mathcal{M}$ and
   has obtained algebraic
properties of Stirling polynomials. He also   proved that  $B_{\mathcal{M}}(m+1)$ is a Hilbert polynomial for any multiset $\mathcal{M}$.
   Analogous problems have been studied for  Legendre-Stirling polynomials \cite{Egge} and  Jacobi-Stirling
polynomials \cite{Gesse2}.    Dzhumadil'daev  and   Yeliussizov  \cite{Dzhumadil'daev}
 provided combinatorial interpretations  of $B_{\mathcal{M}}(m)$ in terms of $P$-partitions, set partitions and permutations.

In   analogy to Stirling permutations, Archer, Gregory, Pennington and Slayden \cite{Archer} introduced   {\em quasi-Stirling} permutations. A permutation $\pi$ of a multiset is said to be a {\em quasi-Stirling
} permutation  if  there does not exist four indices $i<j<k<\ell$ such that $\pi_i=\pi_k$ and $\pi_j=\pi_\ell$. For a multiset $\mathcal{M}$, denote by   $\overline{\mathcal{Q}}_{\mathcal{M}}$ the set of quasi-Stirling permutations of $\mathcal{M}$. For example, if $\mathcal{M}=\{1, 2^2\}$, we have
$
\overline{\mathcal{Q}}_{\mathcal{M}} =\{122, 221, 212\}.
$
 Archer, Gregory, Pennington and Slayden \cite{Archer} noticed that when $\mathcal{M}=\{1^2, 2^2, \ldots, n^2\}$,
 $$
 |\overline{\mathcal{Q}}_{\mathcal{M}}|=n! C_n,
 $$
 where $C_n={1\over n+1}{2n\choose n}$ is the $n$th Catalan number.  They also posed the following intriguing conjecture.
 \begin{conjecture}{\upshape   ( \cite{Archer}) }\label{con1}
  Let $\mathcal{M}=\{1^2, 2^2, \ldots, n^2\}$. The number of $\pi\in \overline{\mathcal{Q}}_{\mathcal{M}}$ with $des(\pi)=n$ is equal to $(n+1)^{n-1}$.
  \end{conjecture}
 Elizalde \cite{Elizalde} confirmed   Conjecture \ref{con1} by establishing a bijection between such quasi-Stirling  permutations  and edge-labeled plane (ordered) rooted trees. He further proved that
 the number of $\pi\in \overline{\mathcal{Q}}_{\mathcal{M}}$ with $des(\pi)=n$ is equal to $((k-1)n+1)^{n-1}$ when $\mathcal{M}=\{1^k, 2^k, \ldots, n^k\}$ and $k\geq 2$.

 The polynomial  $$ \overline{Q}_{\mathcal{M}}(t)=\sum_{\pi\in \overline{\mathcal{Q}}_{\mathcal{M}}}t^{des(\pi)}$$ is called
the  {\em quasi-Stirling polynomial } on the multiset $\mathcal{M}$.
For instance, if $\mathcal{M}=\{1, 2^2, 3\}$, then we have $\overline{Q}_{\mathcal{M}}(t)=t+7t^2+4t^3.$
By employing generating function arguments,  Elizalde \cite{Elizalde}     derived that
  \begin{equation}\label{eq3}
 \sum_{m=0}^{\infty} {m^n\over n+1}{n+m\choose m}t^m={\overline{Q}_{\mathcal{M}}(t)\over (1-t)^{2n+1}},
 \end{equation}
when $\mathcal{M}=\{1^2, 2^2, \ldots, n^2\}$
   and asked for a combinatorial  proof.

    Motivated by  previous results for Stirling polynomials on   multisets,  we    investigate  quasi-Stirling polynomials on multisets      and   derive the following generalization  of   Elizalde's formula (\ref{eq3}).
\begin{theorem}\label{mainth1}
Let $\mathcal{M}=\{1^{k_1}, 2^{k_2}, \ldots, n^{k_n}\}$ and $K=k_1+k_2+\ldots+k_n$ with $k_i\geq 1$. We have
\begin{equation}\label{eqmain}
\sum_{m=0}^{\infty} {m^n\over K-n+1}{K-n+m \choose m} t^m={\overline{Q}_{\mathcal{M}}(t)\over (1-t)^{K+1}}.
\end{equation}
\end{theorem}

  In this paper, we shall provide   a combinatorial proof of (\ref{eqmain}) in terms of certain ordered  labeled trees.  Specializing  $\mathcal{M}=\{1^2, 2^2, \ldots, n^2\}$ in (\ref{eqmain})   implies a combinatorial proof of (\ref{eq3}) in answer to the problem posed by Elizalde \cite{Elizalde}. Notice that  by letting  $\mathcal{M}=\{1,2,\ldots, n\}$, we are led to the well-known  identity \cite{Fota}  involving Eulerian polynomials
  $$
  \sum_{m=0}^{\infty} m^n t^m={A_n(t)\over (1-t)^{n+1}}.
  $$

 Notice that for any multiset $\mathcal{M}=\{1^k, 2,3,\ldots, n\}$ with $k\geq 1$, the quasi-Stirling permutations of $\mathcal{M}$ are equivalent to  Stirling permutations of $\mathcal{M}$.
The following corollary follows immediately from   Theorem \ref{mainth1}.
\begin{corollary}\label{coro1}
Let $\mathcal{M}=\{1^{k_1}, 2^{k_2}, \ldots, n^{k_n}\}$ and $K=k_1+k_2+\ldots+k_n$ with $k_i\geq 1$. Then
    we have $\overline{Q}_{\mathcal{M}}(t)=Q_{\mathcal{M'}}(t)$, where $M'=\{1^{K-n+1}, 2,3,\ldots, n\}$.
    \end{corollary}

It is   well known result of Frobenius that the roots of the Eulerian polynomials
are real, distinct, and nonpositive. B\'{o}na \cite{Bona} proved the real realrootedness  of  Stirling polynomial   $Q_{\mathcal{M}}(t)$ in the case of  $M=\{1^2, 2^2, \ldots, n^2\}$.
    Brenti  \cite{Brenti1}    showed that
the  Striling polynomial   $Q_{\mathcal{M}}(t)$   has only real roots for any multiset $\mathcal{M}$.
Combining  Corollary \ref{coro1} and  Brenti's result,  it is clear  that the quasi-Stirling
polynomial  $\overline{Q}_{\mathcal{M}}(t)$ has  the same  property as  the Stirling polynomial   $Q_{\mathcal{M}}(t)$, which is an extension of Elizalde's result \cite{Elizalde} for $\mathcal{M}=\{1^2, 2^2, \ldots, n^2\}$.
    \begin{corollary}\label{coro3}
    For any multiset $\mathcal{M}$, the quasi-Stirling polynomial $\overline{Q}_{\mathcal{M}}(t)$ has only real roots.
    \end{corollary}
     Unimodal and log-concave sequences arise frequently in combinatorics. It is well known that the coefficients of a   polynomial with nonnegative coefficients and with only real roots are log-concave and that log-concavity implies
unimodality.
    As an application of this
result,  Corollary \ref{coro3} implies the following property of the coefficients of $\overline{Q}_{\mathcal{M}}(t)$.
    \begin{corollary}\label{coro2}
    For any multiset $\mathcal{M}$, the coefficients of $\overline{Q}_{\mathcal{M}}(t)$  are unimodal  and log-concave.
    \end{corollary}
    Notice that specializing $\mathcal{M}=\{1^2, 2^2, \ldots, n^2\}$ in Corollary \ref{coro2}  recovers  Elizalde's result \cite{Elizalde}.

\section{  Proof of Theorem \ref{mainth1}}

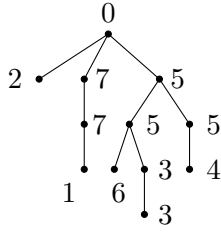
\begin{figure}[H]
\begin{center}
\begin{tikzpicture}[font =\small , scale = 0.4]
\filldraw[fill=black](5.3,1.5)circle(0.1);\filldraw[fill=black](3,0)circle(0.1);\filldraw[fill=black](4.5,0)circle(0.1);\filldraw[fill=black](7,0)circle(0.1);

\filldraw[fill=black](4.5,-1.5)circle(0.1);\filldraw[fill=black](6,-1.5)circle(0.1);\filldraw[fill=black](8,-1.5)circle(0.1);
\filldraw[fill=black](4.5,-3)circle(0.1);\filldraw[fill=black](5.5,-3)circle(0.1);\filldraw[fill=black](6.5,-3)circle(0.1);\filldraw[fill=black](8,-3)circle(0.1);
\filldraw[fill=black](6.5,-4.5)circle(0.1);
\draw(5.3,1.5)--(3,0);\draw(5.3,1.5)--(4.5,0);\draw(5.3,1.5)--(7,0);\draw(4.5,0)--(4.5,-1.5);\draw(4.5,-1.5)--(4.5,-3);\draw(7,0)--(6,-1.5);\draw(7,0)--(8,-1.5);\draw(6,-1.5)--(5.5,-3);
\draw(6,-1.5)--(6.5,-3);\draw(8,-1.5)--(8,-3);\draw(6.5,-3)--(6.5,-4.5);

\coordinate [label=above:$0$] (x) at (5.3,1.55);\coordinate [label=left:$2$] (x) at (2.8,0);\coordinate [label=right:$7$] (x) at (4.5,0);
\coordinate [label=right:$5$] (x) at (7,0);\coordinate [label=right:$7$] (x) at (4.4,-1.5);\coordinate [label=right:$5$] (x) at (6.2,-1.5);\coordinate [label=right:$5$] (x) at (8.2,-1.5);
\coordinate [label=below:$1$] (x) at (4,-3.1);\coordinate [label=right:$3$] (x) at (6.6,-3);\coordinate [label=below:$6$] (x) at (5.65,-3.1);\coordinate [label=right:$4$] (x) at (8.2,-3);
\coordinate [label=right:$3$] (x) at (6.6,-4.5);

\end{tikzpicture}
\end{center}
\caption{ A tree $T\in \mathcal{T}_{
\mathcal{M}}$ with $M=\{1,2,3^{2},4,5^{3},6,7^{2}\}$. }\label{T}
\end{figure}

The objective of this section is to provide a combinatorial proof of Theorem \ref{mainth1}.
Recall that an  {\em ordered} tree is a tree with one designated vertex, which is called the root, and the subtrees of
each vertex are linearly ordered.  In an ordered tree $T$,   the {\em level} of a vertex $v$ in $T$  is defined to be the length of the unique path from the root to $v$.  A vertex $v$ is said to be {\em at odd level} (resp. {\em at even level}) if the level of $v$ is odd (resp. even).
 In order to prove (\ref{eqmain}), we need to consider  a new
class of ordered   labeled trees $T$   which  verify the following properties:
\begin{itemize}
\item[{\upshape (\rmnum{1})}] the vertices are labeled by the elements of the multiset  $\{0\}\cup \mathcal{M}$, where $\mathcal{M}=\{1^{k_1}, 2^{k_2}, \ldots, n^{k_n}\}$ with $k_i\geq 1$;
    \item[{\upshape (\rmnum{2})}] the root is labeled by $0$;
    \item[{\upshape (\rmnum{3})}] for  a vertex  $v$ at odd level,  if $v$ is labeled by $i$, then $v$ has  exactly  $k_i-1$ children  and the children of $v$ have the same label as that of $v$
  in $T$.
\end{itemize}
Let $\mathcal{T}_{\mathcal{M}}$ denote the set of  such ordered    labeled trees. For example,  a tree $T\in \mathcal{T}_{
\mathcal{M}}$ with $\mathcal{M}=\{1,2,3^{2},4,5^{3},6,7^{2}\}$ is illustrated in Figure \ref{T}.

Define the number of cyclic descents of a sequence $\pi=\pi_1\pi_2\ldots\pi_n$ to
be
$$
cdes(\pi)=|\{i\mid \pi_i>\pi_{i+1}\}|
$$
with the convention $\pi_{n+1}=\pi_1$.

In an ordered  labeled tree $T$, let $u$ be a vertex of $T$. Suppose that $u$ has $\ell$ children $v_1, v_2, \ldots, v_{\ell}$ listed from left to right. The number of  cyclic descents of $u$, denoted by $cdes(u)$,  is defined to be  the number  $cdes(uv_1v_2\ldots v_{\ell})$. The number of cyclic descents of $T$  is defined to be
$$
cdes(T)=\sum_{u\in V(T)}cdes(u),
$$
where $V(T)$ denotes the vertex set of $T$.
If $T$ has only one vertex,    let $cdes(T)=1$. For example, if let $T$ be a tree as shown  in Figure \ref{T}, we have $cdes(T)=5$.

In order to prove (\ref{eqmain}), we establish a bijection between $\mathcal{T}_{\mathcal{M}}$ and $
\overline{\mathcal{Q}}_{\mathcal{M}}$.  Let $\pi$ be a permutation, the leftmost (resp. rightmost) entry of $\pi$ is denoted by $first(\pi)$ (resp. $last(\pi)$). Similarly, for an ordered labeled tree $T$, we denote by $first(T)$   (resp. $last(T)$) the leftmost (resp. rightmost) child of the root.
If a tree $T$ has only one vertex, let $first(T)=+\infty$ and $last(T)=-\infty$. Similarly, for an empty permutation $\epsilon$, we   let $first(\epsilon)=+\infty$ and $last(\epsilon)=-\infty$.

\begin{theorem}\label{thphi}
There exists a bijection $\phi$ between $\mathcal{T}_{\mathcal{M}}$ and $
\overline{\mathcal{Q}}_{\mathcal{M}}$  such that $$(cdes, first, last)T= (des, first, last)\phi(T)$$ for any $T\in \mathcal{T}_{\mathcal{M}}$.
\end{theorem}

\pf First we give   a recursive description of   the map $\phi$ from  $\mathcal{T}_{\mathcal{M}}$ to $
\overline{\mathcal{Q}}_{\mathcal{M}}$.
Let $T\in \mathcal{T}_{\mathcal{M}} $. If $T$ has only one vertex, let $\phi(T)=\epsilon$.  Otherwise, suppose that $first(T)=r$.

\noindent Case 1.  The leftmost child of the root is a leaf. Let $T_0$ be the tree obtained from $T$ by removing  the leftmost child of the root together with  the edge incident to it. Define
    $\phi(T)=r\phi(T_0)$.

    \noindent Case 2. The leftmost child of the root has $k$ children. For $1\leq i\leq k$,  let $T_i$ be the subtree rooted at the $i$-th child of the leftmost child of the root.
         Denote by $T'_i$ the tree obtained from $T_i$ by relabeling its root by $0$.
         Let $T_0$ be the tree obtained from $T$ by removing all the subtrees $T_1$, $T_2$, $\ldots,$ $T_k$ and  all the vertices labeled by $r$ together with  the edges incident to them. Define $\phi(T)=r\phi(T'_1)r\phi(T'_2)\ldots r\phi(T'_k)r\ \phi(T_0)$.

It is routine to check that we have $\phi(T)\in \overline{\mathcal{Q}}_{\mathcal{M}} $.  Now we proceed to show that $(cdes, first, last)T= (des, first, last)\phi(T)$  by induction on the cardinality of the multiset $\mathcal{M}$. Obviously, the statement holds when $|\mathcal{M}|=0$. Assume that $(cdes, first, last)T= (des, first, last)\phi(T)$ for any $T\in \mathcal{T}_{\mathcal{M}'}$    with $|\mathcal{M}'|<|\mathcal{M}|$. From the construction of $\phi$, it is easily seen that
$first(T)=first(\phi(T))$.  Moreover, if $T_0$ has only one vertex, one can easily check that  $last(T)=r=last(\phi(T))$.  Otherwise, we have $last(T_0)=last(T)$ and  $last(\phi(T))=last(\phi(T_0))$. By induction hypothesis, we have $last(T)=last(T_0)=last(\phi(T_0))=last(\phi(T))$ as desired.
It remains to show that $cdes(T)=des(\phi(T))$.
 We have two cases.

\noindent Case 1.    The leftmost child of the root is a leaf. In this case,   we have $\phi(T)=r\phi(T_0)$. It is routine to check that
    $$
    cdes(T)=cdes(T_0)+\chi(r>first(T_0)).
    $$
    Here $\chi(S)=1$ if  the statement $S$ is true, and $\chi(S)=0$ otherwise.
    Then by induction hypothesis, we have
    $$
    \begin{array}{lll}
    cdes(T)&=&cdes(T_0)+\chi(r>first(T_0))\\
    &=&des(\phi(T_0))+\chi(r>first(\phi(T_0)))\\
    &=&des(\phi(T)).
    \end{array}
    $$

 \noindent Case 2.  The leftmost child of the root has $k$ children. Then  we have $\phi(T)=r\phi(T'_1)r\phi(T'_2)\ldots r\phi(T'_k)r\ \phi(T_0)$.  It is easily seen that
    $$
    cdes(T)= cdes(T_0)+\chi(r>first(T_0))+\sum_{i=1}^{k}\big(cdes(T'_i)-1+\chi(r>first(T'_i)\big)+\chi(r<last(T'_i)),    $$
     and $des(\phi(T))$ is given by
    $$
    des(\phi(T_0))+\chi(r>first(\phi(T_0)))+\sum_{i=1}^{k}\big(des(\phi(T'_i))-1+\chi(r>first(\phi(T'_i))\big)+\chi(r<last(\phi(T'_i))).
    $$
   By induction hypothesis,  it is easy to verify that $cdes(T)=des(\phi(T))$.

    So far, we have concluded that $\phi$ is map from $\mathcal{T}_{\mathcal{M}}$ to $
\overline{\mathcal{Q}}_{\mathcal{M}}$  satisfying that $$(cdes, first, last)T= (des, first, last)\phi(T)$$ for any $T\in \mathcal{T}_{\mathcal{M}}.$     It is clear that the construction of $\phi$ is reversible, and hence the map $\phi$ is the desired bijection. This completes the proof. \qed

For example, let $T$ be the tree as shown in Figure \ref{T}. By applying the map $\phi$, we get a permutation $\phi(T)=27175633545$.

Relying on Theorem \ref{thphi}, quasi-Stirling polynomials can be reformulated as follows:
 \begin{equation}\label{tree}
 \overline{Q}_{\mathcal{M}}(t)=\sum_{\pi\in \overline{\mathcal{Q}}_{\mathcal{M}}}t^{des(\pi)}=\sum_{T\in \mathcal{T}_{\mathcal{M}}}t^{cdes(T)}.
 \end{equation}

Let $e$ be an edge of an ordered labeled  tree $T$.  If  $e$ connects vertices $u$ and $v$ and $u$ is the parent of $v$, then we call $v$ the {\em south endpoint} of $e$.  An edge whose south endpoint is an unlabeled leaf is called  a {\em half-edge}.

\begin{figure}[h]
\centerline{\includegraphics[width=11cm]{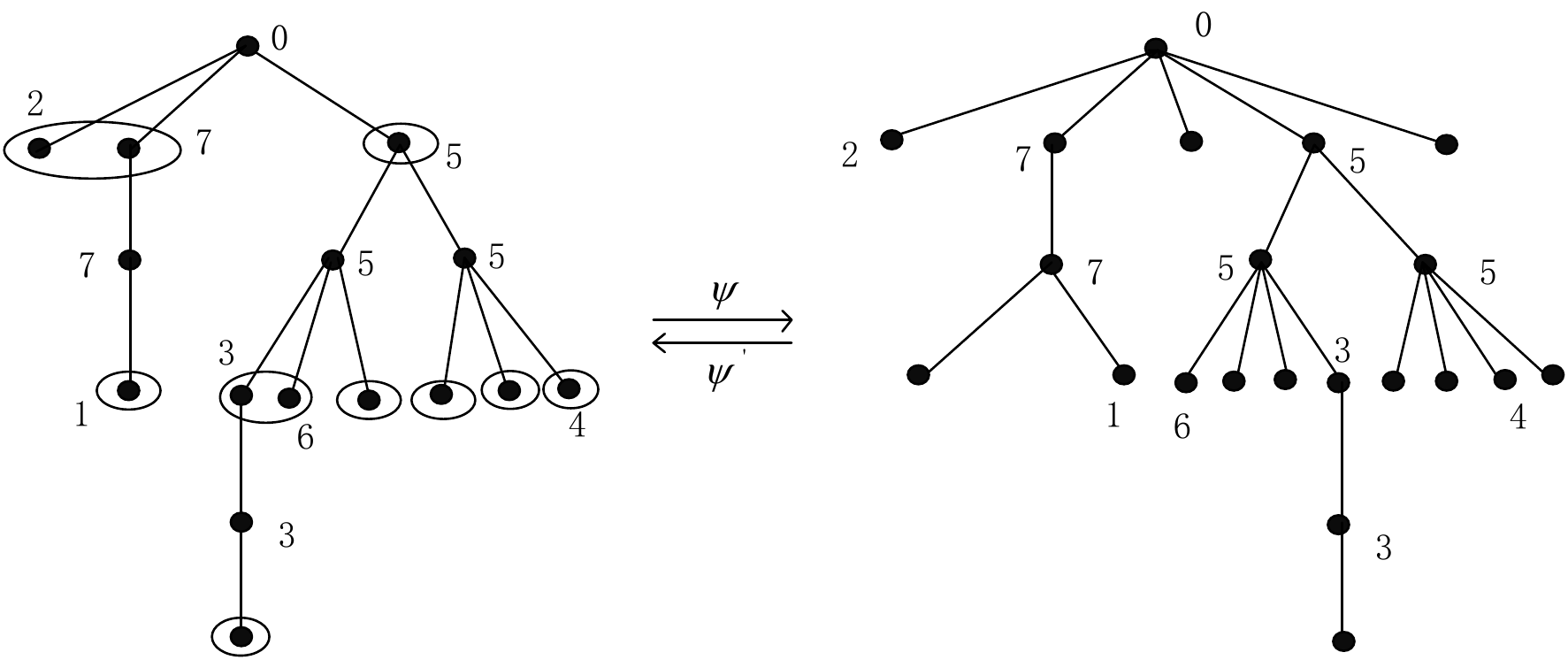}}
\caption{ An example of the bijection $\psi$. }\label{psi}
\end{figure}

Let $T\in \mathcal{T}_{\mathcal{M}}$ and
let $u$ be a vertex of $T$ at even level. Now we may attach some half-edges to the vertex $u$ by obeying the following rules:
\begin{itemize}
\item the inserted half-edges serve as walls to separate the children of $u$ into (possibly empty) compartments such that the vertices in the same compartment are increasing from left to right;
    \item if $v$ is the last child of $u$ and  $v>u$, then there is some half-edge inserted to the right of $v$;
        \item if $v$ is the first  child of $u$ and  $u>v$, then there is some half-edge  inserted to the left  of $v$.
\end{itemize}
   Repeat the above procedure for all the vertices at even level of $T$, we will get a tree with some half-edges inserted. Denote by $\mathcal{T}^*_{\mathcal{M},m}$ the set of  such  trees  with $m$  half-edges    inserted.  For example,  a tree $T\in \mathcal{T}^*_{\mathcal{M},9} $ with $\mathcal{M}=\{1,2, 3^2, 4, 5^3, 6, 7^2\}$ is illustrated  on the right of  Figure \ref{psi}.

   Relying  on $(\ref{tree})$, the coefficient of $t^m$ of the right-hand side of  (\ref{eqmain})  can be interpreted as    the number of ordered   labeled trees in $\mathcal{T}^*_{\mathcal{M},m}$.
 In order to provide a combinatorial proof of (\ref{eqmain}), it suffices to show that
  \begin{equation}\label{eqT*}
 |\mathcal{T}^*_{\mathcal{M},m}|={m^n\over K-n+1}{K-n+m \choose m}
 \end{equation}
   for $\mathcal{M}=\{1^{k_1}, 2^{k_2}, \ldots, n^{k_n}\}$ and $k_1+k_2+\ldots+k_n=K$ with $k_i\geq 1$.
To this end, we need to consider a class of ordered   labeled trees satisfying that
\begin{itemize}
\item[{\upshape (\rmnum{1})}] the vertices are labeled by the elements of the multiset  $\{0\}\cup \mathcal{M}$ and a leaf at odd level  may  be unlabeled, where $\mathcal{M}=\{1^{k_1}, 2^{k_2}, \ldots, n^{k_n}\}$;
    \item[{\upshape (\rmnum{2})}] the root is labeled by $0$;
    \item[{\upshape (\rmnum{3})}] for each  vertex  $v$ at odd level,  if $v$ is labeled by $i$, then $v$ has  exactly  $k_i-1$ children  and the children of $v$ have the same label as that of $v$;
\item[{\upshape (\rmnum{4})}] for each internal vertex  $v$ at even level, the children of $v$  are separated  into blocks  satisfying  that each block  is  either     an unlabeled leaf or  consists of    children which are  increasing from left to right in $T$.
\end{itemize}
  Denote by $\mathcal{BT}_{\mathcal{M},m}$ the set of such ordered    labeled trees      with $m$ blocks.     For example,  a tree $T\in \mathcal{BT}_{\mathcal{M},9} $ with $\mathcal{M}=\{1,2, 3^2, 4, 5^3, 6, 7^2\}$ is illustrated  on the left of  Figure \ref{psi}, where the vertices  in the same block are grouped together by a circle.

 The bijective proof of  (\ref{eqT*}) relies on the following two results.
 \begin{theorem}\label{thpsi}
 There is a bijection $\psi:\mathcal{BT}_{\mathcal{M},m}\rightarrow \mathcal{T}^{*}_{\mathcal{M},m}$.
\end{theorem}
 \begin{theorem}\label{thBT}
 Let $\mathcal{M}=\{1^{k_1}, 2^{k_2}, \ldots, n^{k_n}\}$ and $k_1+k_2+\ldots+k_n=K$ with $k_i\geq 1$. We have
  \begin{equation}\label{eqBT}
 |\mathcal{BT}_{\mathcal{M}, m}|= {m^n\over K-n+1}{K-n+m\choose m}.
 \end{equation}
 \end{theorem}

 \noindent{\bf Proof of Theorem \ref{thpsi}.} Let $T\in \mathcal{BT}_{\mathcal{M},m}$ and let $u$ be an internal vertex at even level of $T$. A block of $u$ is said to be  {\em trivial} if the block consists of an unlabeled leaf.
The map $\psi: \mathcal{BT}_{\mathcal{M},m} \rightarrow  \mathcal{T}^{*}_{\mathcal{M},m}$ can be described as follows.   For an  internal  vertex $u$ at even level of $T$, attach a half-edge to $u$   to the right of   the rightmost vertex of  each non-trivial block  of $u$. Assume that  the vertices of the leftmost block of $u$  are  $v_1, v_2, \ldots, v_{\ell}$  with $v_1<v_2<\ldots<v_t<u<v_{t+1}<\ldots<v_{\ell}$. Then
  rearrange   $v_1$, $v_2$, $\ldots$, $v_t$  as the last $t$ children of $u$, which are arranged in  increasing order from left to right.  Repeat the above procedure for each internal vertex at even level  of $T$. Define $\psi(T)$ to be the resulting tree. It is routine to check that $\psi(T)\in \mathcal{T}^{*}_{\mathcal{M},m}$ and thus $\psi$ is well defined.

  In order to show that $\psi$ is bijection, we proceed to describe a map $\psi':\mathcal{T}^{*}_{\mathcal{M},m} \rightarrow  \mathcal{BT}_{\mathcal{M},m}$. Let   $T\in \mathcal{T}^{*}_{\mathcal{M},m}$ and let $u$ be an internal vertex at even level of $T$. Assume that the  children of $u$ to the right of  the last half-edge incident to $u$ are  $v_1, v_2, \ldots, v_t$.  Clearly, we have $v_1<v_2<\ldots<v_t<u$. Rearrange the vertices    $v_1, v_2, \ldots, v_t$  as   the leftmost $t$ children of $u$, arranged in increasing order from left to right.  Then
   the half-edges serve as walls to separate the children of $u$ into (possibly empty) compartments such that the vertices in the same compartment are increasing from left to right.
      Remove each half-edge $e$ incident to $u$   if  the compartment   to the left of  $e$ is non-empty. Then let   each non-empty compartment  (resp.  each unlabeled leaf) be a block of $u$.    Repeat the above procedure for each internal vertex at even level  of $T$. Define $\psi'(T)$ to be the resulting tree.
  It is easy to check that $\psi'(T)\in \mathcal{BT}_{\mathcal{M},m} $. Moreover, the maps $\psi$ and $\psi'$ are inverses of each other, and thus $\psi$ is a bijection as desired. This completes the proof. \qed

For example,  let  $T$ be a tree in $\mathcal{BT}_{\mathcal{M}, 9} $ with $M=\{1,2,3^{2},4,5^{3},6,7^{2}\}$   as illustrated  on the left of   Figure \ref{psi}. By applying the map $\psi$, we  obtain a tree $\psi(T)\in \mathcal{T}^*_{\mathcal{M}, 9} $ as shown on the   right of   Figure \ref{psi}.

In the following, we shall  present a bijective enumeration of trees in $\mathcal{BT}_{\mathcal{M},m}$ by exhibiting   a Pr$\ddot{u}$fer-like code.

Let $\mathcal{M}=\{1^{k_1}, 2^{k_2}, \ldots, n^{k_n}\}$ and  $I=\{i\mid k_i>1\}$.  For $k_i>1$, let
$\mathcal{M}_i=\{i_1, i_2, \ldots, i_{k_i-1}\}$.
Denote by $\mathcal{P}_{\mathcal{M},m}$ the set of pairs $(P, S)$ where
\begin{itemize}
\item $P$ is an $m$-multiset on the  set $\{0\} \bigcup(\bigcup_{i\in I}\mathcal{M}_i)$;
\item $S$ is a sequence $(a_1, b_1)(a_2, b_2) \ldots (a_n, b_n)$  satisfying  that  $a_i \in P$  and  $b_i$ is a positive integer not larger than the number of occurrences  of $a_i$ in $P$;
    \item $a_n=0$
\end{itemize}

We need a combinatorial identity before we get the enumeration of the pairs in $ \mathcal{P}_{\mathcal{M},m}$.
\begin{lemma}\label{lemiden}
 For  $m\geq 1$, we have
 \begin{equation}\label{eqiden}
 \sum_{\ell=1}^{m}\ell{n+m-\ell-1\choose n-1}={n+m\choose n+1}.
\end{equation}
\end{lemma}
\pf  It is apparent that  the statement  holds for $m=1$. Now we assume $m\geq 2$. The right-hand side of (\ref{eqiden}) can be interpreted as the number of $(n+1)$-subsets $P$ of $[n+m]=\{1,2, \ldots, n+m\}$.  Suppose that the second largest number appear in $P$ is equal to $\ell+1$ for some $1\leq \ell\leq m$. Then  we have $\ell$ ways to choose the minimal element of $P$ from $[\ell]$, and we can choose the remaining $n-1$ elements of $P$ from $[n+m]\setminus [\ell+1]$ in  ${n+m-\ell-1\choose n-1}$ ways. Therefore,   the total number of ways to choose such a subset $P$ is given by $\ell{n+m-\ell-1\choose n-1}$. Then the result follows, completing the proof. \qed

\begin{lemma}\label{lemP}
Let $\mathcal{M}=\{1^{k_1}, 2^{k_2}, \ldots, n^{k_n}\}$ and $K=k_1+k_2+\ldots+k_n$ with $k_i\geq 1$. We have
$$
|\mathcal{P}_{\mathcal{M},m}|={m^n\over K-n+1}{K-n+m\choose m}.
$$
\end{lemma}
  \pf  Given a positive integer $\ell$ with $\ell \leq m$, we can generate a pair $(P,S)\in \mathcal{P}_{\mathcal{M},m}$ in which the number of occurrences of $0$ in $P$ equals $\ell$. It is apparent that we have ${K-n+m-\ell-1\choose m-\ell}$ ways to choose such a multiset $P$. Once $P$ is determined, we can choose $b_n$ in  $\ell$ ways, and  each $(a_i, b_i)$ in  $m$ ways for $1\leq i\leq n-1$. Then the total number of different ways to choose a   pair   $(P,S)\in \mathcal{P}_{\mathcal{M},m}$ in which the number of occurrences of $0$ in $P$ equals $\ell$ is given by $m^{n-1}\ell{K-n+m-\ell-1\choose m-\ell}$. This yields that
  $$
  |\mathcal{P}_{\mathcal{M},m}|=m^{n-1}\sum_{\ell=1}^{m}\ell{K-n+m-\ell-1\choose m-\ell}.
  $$
By  (\ref{eqiden}), we have
 $$
 \begin{array}{lll}
  |\mathcal{P}_{\mathcal{M},m}|&=&m^{n-1}\sum_{\ell=1}^{m}\ell{K-n+m-\ell-1\choose m-\ell}\\
  &=& m^{n-1}{K-n+m\choose K-n+1}\\
  &=&{ m^{n}\over K-n+1}{K-n+m\choose m}\\
  \end{array}
  $$
  as desired, completing the proof. \qed

For a tree $T\in \mathcal{BT}_{\mathcal{M},m} $, a vertex $v$ is said to be {\em even} (resp. {\em odd} ) if $v$ is at even (resp. odd) level.
\begin{lemma}\label{lemtheta}
There exists a bijection  $\theta:\mathcal{BT}_{\mathcal{M},m} \rightarrow \mathcal{P}_{\mathcal{M},m}$.
\end{lemma}
\pf Let $T\in \mathcal{BT}_{\mathcal{M},m} $.
First we obtain a     tree $T'$ on the vertex set $\mathcal{M}\cup \{0\}$ by   removing  all the half-edges from $T$.      Starting from $T_0=T'$, for $i=1,2,\ldots, n$, let $T_i$ be the tree obtained from $T_{i-1}$ by deleting the largest odd vertex $l_i$  together with  its children  if  all the children of the odd vertex $l_i$ are leaves  or the odd vertex $l_i$ is a leaf.     Set $\theta(T)=(P, S)$ where
\begin{itemize}
\item $P$ is a mutiset on the set $\{0\} \bigcup(\bigcup_{i\in I}\mathcal{M}_i)$  such that the
\begin{itemize}
\item for all $x\in [n]$, the number of occurrences of $x_j$ in $P$ equals  $c$ if and only if the  children of  the $j$-th child of the odd vertex $x$  are separated  into  $c$ blocks  (both children and blocks are ordered from left to right);
    \item the number of occurrences of $0$ in $P$ equals  $c$ if and only if the  children of  the root $0$  are separated  into  $c$ blocks.
\end{itemize}
\item $S$ is a sequence $(a_1, b_1)(a_2, b_2) \ldots (a_n, b_n)$  satisfying  that
 \begin{itemize}
  \item $a_i=x_j$ for some $x\in [n]$  if and only if  the odd vertex $l_i$ lies in the $b_i$-th block of the $j$-th child of the odd vertex $x$;
      \item $a_i=0$  if and only if  the vertex $0$ is the parent of the odd vertex $l_i$ and the odd vertex  $l_i$ lies in the $b_i$-th block of   the vertex $0$ in $T$.
 \end{itemize}
\end{itemize}
It is apparent that we have $\theta(T)\in \mathcal{P}_{\mathcal{M},m}$.

Conversely, given a pair $(P, S)\in \mathcal{P}_{\mathcal{M},m}$, we can recover a tree $T\in \mathcal{BT}_{\mathcal{M},m} $  by the following procedure.
\begin{itemize}
\item[(a)] For $1\leq i\leq n$, let $l_i$ be the largest element of $[n]\cup\{0\}\setminus \{l_1, l_2, \ldots, l_{i-1}\}$ such that $(l_i)_j$ does not appear in the sequence $(a_{i}, a_{i+1}, \ldots, a_{n})$ for all $1\leq j\leq k_{l_i}-1$;
    \item[(b)] If $a_i=x_j$ for some $x\in [n]$, then let   $l_i$ be an odd  vertex lying in the $b_i$-th block of    the $j$-th child of  the odd vertex $x$.  Let the odd vertex $l_i$ have  $k_{l_i}-1$ children  and let all the children of the odd vertex $l_i$ be labeled by $l_i$.
         \item[(c)] If $a_i=0$, then  let $l_i$ be an odd vertex lying in the $b_i$-th block of    the vertex $0$.  Let the odd vertex  $l_i$ have   $k_{l_i}-1$ children  and let all the children of the odd vertex $l_i$ be labeled by $l_i$.
    \item[(d)] Let $(x,y)\neq (a_i, b_i)$ for any $i\in [n] $,  where  $x\in P$ and  $y$ is a positive integer not larger than the number of the occurrences  of $x$ in $P$.
           \begin{itemize}
           \item  If $x=i_j$ for some $1\leq i\leq n$ and $1\leq j\leq k_i-1$, we  attach a half-edge to the $j$-th child of the odd vertex $i$   such that the new inserted   unlabeled  leaf lies in its $y$-th block;
                \item  If $x=0$, we attach a  half-edge to  the vertex  $0$ such that the new inserted   unlabeled  leaf lies in its $y$-th block.
           \end{itemize}
\end{itemize}
 It is clear that the map $\theta$ is reversible, and therefore  the map $\theta$ is a bijection. This completes the proof.  \qed

\begin{figure}[h]
\centerline{\includegraphics[width=15cm]{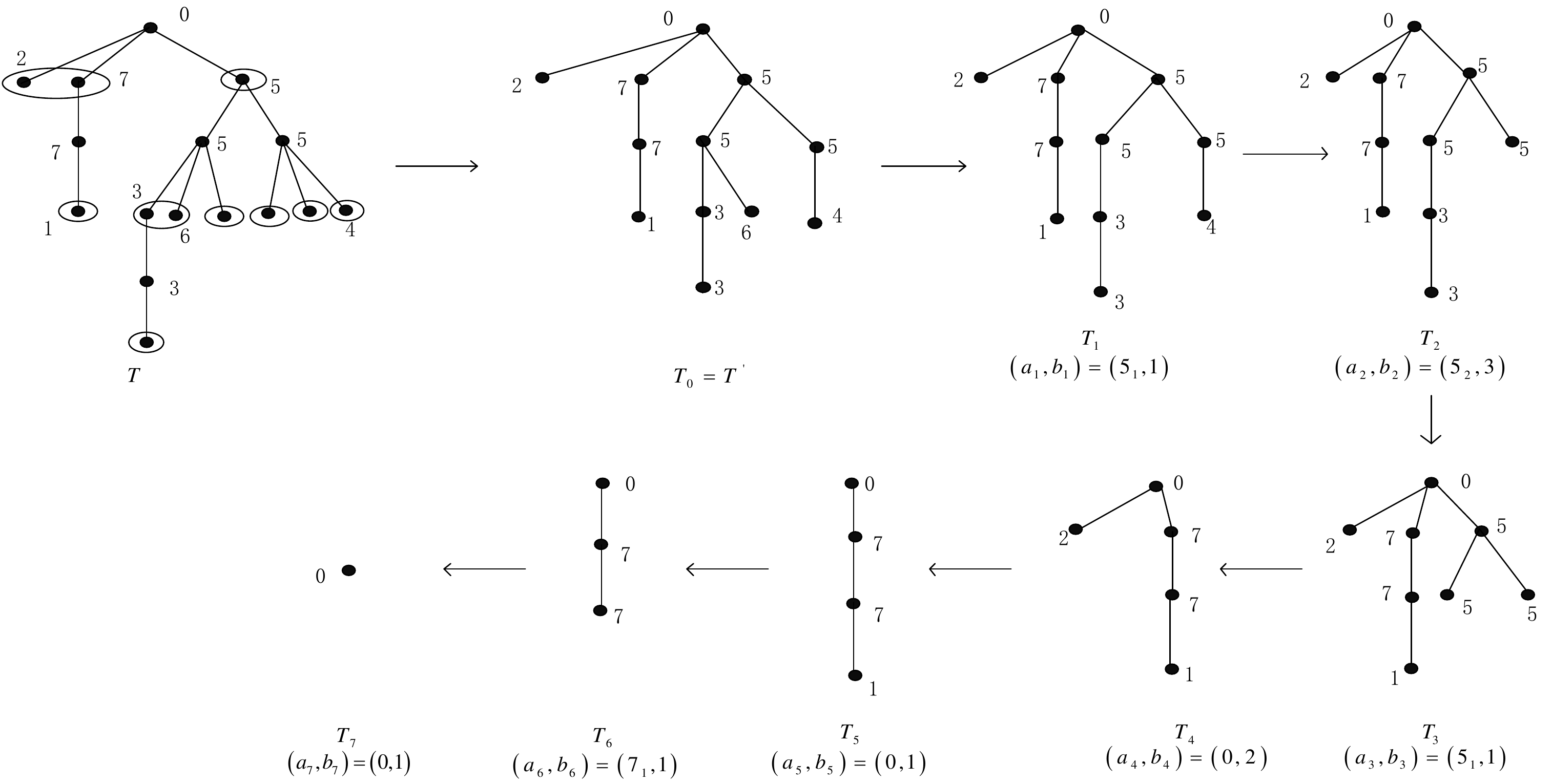}}
\caption{ An example of the bijection $\theta$. }\label{theta}

\end{figure}

For example, let $T$ be a tree in $\mathcal{BT}_{\mathcal{M},9}$ with $ \mathcal{M}=\{1,2, 3^2, 4, 5^3, 6,7^2\}$ as shown in Figure \ref{theta}. By applying the map $\theta$, we get a pair $\theta(T)=(P, S)\in \mathcal{P}_{\mathcal{M},9}$, where  $$P=\{0^2, 3_1, 5_1^2, 5_2^3, 7_1\}$$ and
$$
S=(5_{1},1)(5_{2},3)(5_{1},1)(0,2)(0,1)(7_{1},1)(0,1).
$$

Combining Lemmas \ref{lemP} and \ref{lemtheta}, we are led to a combinatorial  proof
 Theorem \ref{thBT}.  Combining  Theorems \ref{thpsi} and \ref{thBT}, we obtain a combinatorial proof of (\ref{eqT*}), thereby providing a combinatorial proof of Theorem \ref{mainth1}.

\vspace{0.5cm}
 \noindent{\bf Acknowledgments.}
 This work was supported by  the National Natural Science Foundation of China (12071440).


\end{document}